\documentclass[12pt]{amsart}
\usepackage{etoolbox} % to patch the style for subsubsections
\patchcmd{\subsubsection}{\normalfont}{\bfseries\textnormal}{}{} %(aleL) subsections as sections; boldfaced font

\usepackage[vmargin=2cm,hmargin=2cm]{geometry} 

\usepackage{enumitem}
\usepackage{caption} 
\usepackage{array}
\usepackage{ifthen}
\usepackage{url}
\usepackage{amssymb,amsmath, amsfonts}
\newif\iffull
\usepackage[x11names]{xcolor}
\usepackage{graphicx}

\newtheorem{theorem}{Theorem}[section]

\newtheorem{lemma}{Lemma}[section]

\newtheorem{question}{Question}[section]
\theoremstyle{definition}
\newtheorem{definition}{Definition}[section]
\newtheorem{example}{Example}[section]

\newtheorem{proposition}{Proposition}[section]

\newtheorem{remark}{Remark}[section]

\numberwithin{equation}{section}

\def\Q{\mathbb Q}
\def\Z{\mathbb Z}

\newcommand{\klammern}[4][]%
{\ifthenelse{\equal{#1}{}}{\left#2}{\csname#1\endcsname#2}%
#4\ifthenelse{\equal{#1}{}}{\right#3}{\csname#1\endcsname#3}}

\date{\today}
\newcommand{\eps}{\varepsilon}
\usepackage{comment}
\usepackage{hyperref}
\hypersetup{hidelinks}
\usepackage[english]{babel}

\newcommand{\genneighbor}{$10^4$}
\makeatletter
\renewcommand*\env@cases[1][1.2]{%
  \let\@ifnextchar\new@ifnextchar
  \left\lbrace
  \def\arraystretch{#1}%
  \array{@{}l@{\quad}l@{}}%
}
\makeatother

\newenvironment{Proof}[1][Proof]{\par\noindent\textbf{#1.}~}
{\hfill\qedsymbol\smallskip\par}
\begin{document}

\title[Sequences of
integers generated by two fixed primes]{Sequences of integers generated
by\\ two fixed primes}

\author{Alessandro Languasco}
\address{Dipartimento di Matematica ``Tullio Levi-Civita'', 
Universit\`a di Padova, Via Trieste, 63, 35121 Padova, Italy}
\email{alessandro.languasco@unipd.it}

\author{Florian Luca}
\address{School of Mathematics, University of the Witwatersrand, P.\,O.\,Box Wits 2050, South Africa}
\email{florian.luca@wits.ac.za}

\author{Pieter Moree}
\address{Max-Planck-Institut f\"ur Mathematik, Vivatsgasse 7, D-53111 Bonn, Germany}
\email{moree@mpim-bonn.mpg.de}

\author{Alain Togb\'e}
\address{Department of Mathematics and Statistics, Purdue University Northwest, 2200 169th Street, Hammond, IN 46323, USA}
\email{atogbe@pnw.edu}

%    General info
%\date{Received by the editor March 3, 1998 and, in revised form, April 28, 1998.}


\maketitle
\dedicatory{\centerline{\textit{In honorem Robert Tijdeman annos LXXX nati}}}

\begin{abstract}
Let $p$ and $q$ be two distinct fixed prime numbers and $(n_i)_{i\geq 0}$ the sequence of 
consecutive integers of the form $p^a\cdot q^b$ with $a,b\ge 0$. Tijdeman gave a lower bound (1973)
and an upper bound (1974) for the 
gap size $n_{i+1}-n_i$, with each bound containing  an unspecified exponent and 
implicit constant. We will explicitly bound these four quantities.
Earlier Langevin (1976) gave weaker estimates for (only) the exponents.
\par Given a real number $\alpha>1$, there exists a smallest number $m$ such that for every 
$n\ge m$, there exists
an integer $n_i$ in 
$[n,n\alpha)$. 
Our effective version of Tijdeman's result  immediately implies an upper bound for $m$, 
which using the Koksma-Erd\H{o}s-Turan inequality we will improve on.
We present a fast algorithm to determine $m$  when $\max\{p,q\}$ is not too large 
and demonstrate it with numerical material. In an appendix we 
explain, given $n_i$, how to efficiently determine both $n_{i-1}$ and $n_{i+1}$, something closely related to work of 
B\'erczes, Dujella and Hajdu. 
\end{abstract}

%&&&&&&&&&&&&&&&&&&&&&&&&&%
\section{Introduction}\label{sec1}
%&&&&&&&&&&&&&&&&&&&&&&&&&%
Given a set $S=\{p_1,\ldots,p_s\}$ of 
primes, the numbers
$p_1^{e_1}\cdots p_s^{e_s}$ with non-negative exponents are called \emph{$S$-units}.
There are many Diophantine equations involving $S$-units, see for
example Evertse et al.\,\cite{EGST}.
In this note we 
are interested in the distribution of $S$-units.
\par In 
case $p_1,\ldots,p_s$ are the first $s$ prime numbers, the counting function of $S$-units up to $x$ is denoted by $\psi(x,p_s)$, and the $S$-units are called $p_s$--friable or $p_s$--smooth. There is an extensive literature on these numbers, see
Hildebrand and Tenenbaum \cite{HT} for a nice survey. The gaps between
consecutive friable numbers were studied by Heath-Brown \cite{Roger}.
Tijdeman and Meijer \cite{TM} studied the 
distribution of the greatest common divisor of 
consecutive $S$-units (for general $S$).
We will concentrate solely on the case where 
$S$ has two elements.
\begin{definition}
\label{def1}
Let $p,q$ be two primes with 
$p<q$.
We let $(n_i)_{i\geq 0}$ be the sequence of 
consecutive integers of the form $n=p^a\cdot q^b$ with $a,b\ge 0$.
\end{definition}
Put $N_{p,q}(x):=\#\{n_i\le x\}$. Note that $N_{p,q}(x)=\#\{(e,f)\in \mathbb Z_{\ge 0}^2: e\log p +f\log q\le \log x\}$, that is it equals the number of lattice points inside a 
rectangular triangle with sides of length 
$\log x/\log p$, 
respectively $\log x/\log q$. The area of this
triangle is
\begin{equation}
\label{asymptotic}    
\frac{(\log x)^2}{2\log p\cdot \log q}.
\end{equation}
It is known that 
$$N_{p,q}(x)=\frac{
\log px\cdot \log qx}{2\log p\cdot \log q}+o(\log x),$$
where the error term is worse if $\log px\cdot \log qx$ is replaced by $(\log x)^2$. Interestingly, Ramanujan in his famous first letter to Hardy (Jan.\,16th, 1913), claims (in case $p=2$ and $q=3$) the latter main term as an approximation, rather than the trivial \eqref{asymptotic}. For more details and references the reader is referred to Moree \cite{fromto}.

B\'erczes, Dujella, and Hajdu \cite{BDH-2014}, for any term $n_i$ determined $n_{i+1}$,  
at least in principle, without enumerating all terms of the 
sequence, and they 
gave an efficient algorithm to 
find $n_{i+1}$ explicitly. They do so by analyzing the behavior of the continued fractions of $\log p/ \log q$ (see 
Theorems 2.2 and 2.3 of \cite{BDH-2014}). 
In the appendix we present a shorter reproof and in addition show how to also efficiently find $n_{i-1}$ 
explicitly.

As early as 1908, Thue \cite{Thue} gave an ineffective
proof that the gap size $n_{i+1}-n_i$ tends to infinity, which was made
effective by Cassels \cite{Cassels} in 1960.
Tijdeman \cite{Tijdeman-1,Tijdeman-2}, half a 
century ago, gave bounds for the gap size. 
He derived these by making use of estimates of Baker \cite{Baker} for linear forms in logarithms of algebraic numbers, which had become available a few years earlier.
\begin{theorem}\label{le:Tijdeman}
Let $(n_i)_{i\geq 0}$ be as in 
Definition \ref{def1}. There exist effective constants $C_1$ and $C_2$ such that
$$
\frac{n_i}{(\log n_i)^{C_1}} \ll_{p,q}  n_{i+1}-n_i \ll_{p,q}  \frac{n_i}{(\log n_i)^{C_2}}.
$$
The constants $C_1,C_2$ and the two implicit constants all may depend on $p$ and $q$.
\end{theorem}

\par  Langevin \cite{Langevin}, 
soon after Tijdeman's work appeared, gave 
effective bounds for the constants $C_1$ and $C_2$.  
\begin{theorem}\label{le:Langevin}
Under the conditions of Theorem \ref{le:Tijdeman}, we have
$$
(2^{126} \log p\cdot \log q\cdot \log \log p)^{-1}  < C_2 < C_1 < 2^{126}\log p\cdot \log q\cdot \log \log p.
$$
\end{theorem} 

We will improve on this result and also give explicit bounds for the implicit constants in Theorem 
\ref{le:Tijdeman}.
\begin{theorem}
\label{maintijdeman}
Assuming $n_i\ge 3$, we have 
$$
C_3\frac{n_i}{(\log n_i)^{C_1}}<n_{i+1}-n_i<C_4\frac{n_i}{(\log n_i)^{C_2}},
$$ 
where $C_1=2\cdot 10^9\log p\cdot \log q$, $C_2=C_1^{-1}$, $C_3=(\log p)^{C_1}$, $C_4=8q$.
\end{theorem}

Langevin used a result of van der Poorten \cite{Poorten} to obtain 
Theorem \ref{le:Langevin}. 
In Section \ref{sec2} we will use a celebrated result of 
Matveev (Theorem \ref{Matveev11}) to prove
Theorem \ref{maintijdeman}. This helps to 
decrease the 
coefficient of the bound $C_2$ from a $38$ 
digit number 
to a $9$ digit one.
\par Very recently Stewart \cite{Stewart} derived an analogue 
of Theorem \ref{le:Tijdeman} for increasing sequences $(n_i)$ such that the largest prime
factor of $n_i$ is at most $y(n_i)$, with $y(x)$ a non-decreasing
function slowly tending to infinity.
\par In Section \ref{sec3}, we consider Bertrand's Postulate 
type results for our sequence. 
\begin{definition}
\label{def:xpalpha}
Given any real number $\alpha>1$, $n_{p,q}(\alpha)$
is the smallest integer
such that every interval 
$[n,n\alpha)$ with $n\ge n_{p,q}(\alpha)$  
contains an integer of the form 
$p^a\cdot q^b$ for every integer $n\ge n_{p,q}(\alpha)$.
\end{definition}
Bertrand's Postulate type results arise on taking $\alpha=2$. 
Theorem \ref{maintijdeman} gives right away that if $\alpha\in (1,p]$ and 
\begin{equation}
\label{directbound}    
n_i>\exp\Bigl(\Bigl(\frac{4q}{\alpha-1}\Bigr)^{2\cdot 10^9\log p\cdot \log q}\Bigr),
\end{equation}
then $n_{i+1}\in (n_i,\alpha n_i)$.
In Section \ref{sec3} we will (slightly) improve on this lower bound and show that this conclusion already holds if
$$n_i>\exp\Bigl(2\log p\Bigl(\frac{60\log q}{\log \alpha}\Bigr)^{10^9\log p\cdot \log q}\Bigr).$$
In Section \ref{sec4} we present an algorithm for efficiently determining $n_{p,q}(\alpha)$ and present some of its outputs. In particular, 
in Section \ref{sec:lucasdisc} we use our algorithm to advance the understanding of the so-called \emph{discriminator} of an infinite family of second-order recurrent sequences first studied by Faye, Luca and Moree \cite{FLM} and more
recently by Ferrari, Luca 
and Moree \cite{FerLM}.

%&&&&&&&&&&&&&&&&&&&&&&&&&&&&&&&&&&&&&%
\section{Proof of Theorem \ref{maintijdeman}}\label{sec2}
%&&&&&&&&&&&&&&&&&&&&&&&&&&&&&&&&&&&&&
We need linear forms in logarithms. For any non-zero algebraic number $\eta$ of degree $d$ over $\Q$, whose minimal polynomial
over $\Z$ is $a\prod_{i=1}^d \bigl(X-\eta^{(i)} \bigr)$ (with $a>0$), we denote by
$$
h(\eta) = \frac{1}{d} \Bigl( \log a + \sum_{i=1}^d \log\max\bigl(1,  \vert \eta^{(i)} \vert  \bigr)\Bigr)
$$
the usual \emph{absolute logarithmic height} of $\eta$. 
If $\eta_1$ and $\eta_2$ are algebraic numbers, then we have the basic properties
\begin{align*}
 h(\eta_1 \pm \eta_2) &\leq h(\eta_1)+ h(\eta_2) +\log2,\\
h(\eta_1\eta_2^{\pm}) &\leq h(\eta_1) + h(\eta_2),\\
h(\eta_1^j)&= \vert j \vert h(\eta_1),
\end{align*} 
where $j$ is any integer.
\par We recall Matveev's main theorem \cite{Mat} in a version 
due to Bugeaud, Mignotte and Siksek \cite[Thm.\,9.4]{BMS:2006}. It applies to algebraic numbers, but we recall it here only for rational numbers. 
For a rational number $\gamma=r/s$ with coprime integers $r$ and $s>0$,
let $h(\gamma):=\max\{\log  \vert r \vert , \log s\}$ be its \emph{naive height}.

\begin{theorem}[Matveev's theorem]
\label{Matveev11} Let $\gamma_1,\ldots,\gamma_k$ be positive rational
numbers, let $b_1,\ldots,b_k$ be non-zero integers, and assume that
\begin{equation}
\label{eq:Lambda}
\Lambda:=\gamma_1^{b_1}\cdots\gamma_k^{b_k} - 1,
\end{equation}
is non-zero. 
For every 
real number $B\geq\max\{ \vert b_1 \vert , \ldots,  \vert b_k \vert \}$ we have
$$
\log  \vert \Lambda \vert  > -1.4\cdot 30^{k+3}\cdot k^{4.5}(1+\log B)\,
h(\gamma_1)\cdots h(\gamma_k).
$$
\end{theorem}

Note that $n_1=1,~n_2=p$. Let $n_i=p^{u} \cdot q^v$. We assume that $n_i\ge 3$ (this holds for $i=2$ in all cases except if $p=2$, in which case we assume that $i\ge 3$). The lower bound on $n_{i+1}-n_i$ follows right--away from Theorem \ref{Matveev11}. Indeed, let 
$n_{i+1}=p^{u'} \cdot q^{v'}$. Note that $\max\{u,v,u',v'\}\le 2\max\{u,v\}$ (since certainly $n_i^2>n_i$ is a $\{p,q\}$-unit, clearly $n_i^2\ge n_{i+1}$). Then 
$$
n_{i+1}-n_i=n_i\bigl(p^{u'-u}q^{v'-v}-1\bigr).
$$
By Theorem \ref{Matveev11}
the right--hand side is
bounded below by 
$$
n_i\,\exp\bigl(-1.4\cdot 30^5\cdot 2^{4.5} (1+\log (2\max\{u,v\}))\log p \cdot \log q\bigl).
$$
Since $\max\{u,v\}\ge 1$ and 
for $x\ge 2$ we have $1+\log x\le 2.5\log x$, we obtain
$$
n_{i+1}-n_i> \frac{n_i}{(2\max\{u,v\})^{2\cdot 10^9\log p\cdot \log q}}.
$$
Since $n_i=p^{u} \cdot q^v\le (pq)^{\max\{u,v\}}$, we get that $\max\{u,v\}\ge \log n_i/\log (pq)$. Thus, 
$$
n_{i+1}-n_i>C_3\frac{n_i}{(\log n_i)^{C_1}},
$$
where 
$$
C_1=2\cdot 10^9\log p\cdot \log q,\qquad C_3=(0.5\log pq)^{C_1}>(\log p)^{C_1}.
$$
In order to obtain an upper bound for $n_{i+1}-n_i$ we proceed as in Langevin \cite{Langevin} and
make use of the sequence
$(r_k/s_k)_{k\ge 0}$ 
of \emph{convergents} of
$$\theta_{p,q}:= \frac{\log p}{\log q}.$$
Suppose  first that $p^{u}>q^v$. We assume that $\ell$ is the index verifying $s_{\ell}\le u<s_{\ell+1}$. Since $p<q$, it follows that $a_0=0$ and $a_1=\lfloor \log q/\log p\rfloor$. So, if $\ell=0$, then we have $u<\log q/\log p$, so $p^{u}<q$. Thus, $v=0$, $n_{i+1}\le q$ and the inequalities hold provided that the multiplicative constant $C_4$ implied by the right Vinogradov symbol $\ll $ is  taken to be at least $q(\log q)^{C_2}$.

Assume next that $\ell\ge 1$. 
One of the rational numbers $r_{\ell}/s_{\ell}$ and $r_{\ell+1}/s_{\ell+1}$ is at least 
$\theta_{p,q}$.
Choose $\ell'\in \{\ell,\ell+1\}$  such that $r_{\ell'}/s_{\ell'}>\theta_{p,q}$. By construction 
$$
p^{u-s_{\ell'}} \cdot q^{v+r_{\ell'}}
$$
is integer larger than $p^{u} \cdot q^v$, and 
so it is at least $n_{i+1}$. 
Thus, we obtain
\begin{equation}
\label{starting}
\log\Bigl(\frac{n_{i+1}}{n_i}\Bigr)
\le 
\vert \Lambda_{\ell'} \vert,
\end{equation}
where 
$$\Lambda_k:=  s_k\log p-r_k\log q.$$
By  \cite[Thm.\,171]{HW} we have
\begin{equation}
\label{cfrac-approx}
\vert  \Lambda_{k}\vert <\frac{\log q}{s_{k+1}}.
\end{equation}
Combining \eqref{starting} and \eqref{cfrac-approx} with $k=\ell' = \ell$ gives
\begin{equation}
\label{eq:Case1}
\log\Bigl(\frac{n_{i+1}}{n_i}\Bigr)<\frac{\log q}{s_{\ell+1}}.
\end{equation}
Next, assume $k=\ell'=\ell+1$. Now 
\eqref{cfrac-approx} gives 
\begin{equation}
\label{eq:2}
%\frac{s_{\ell+1}}{\log q}<\frac{1}{|\Lambda|},
\frac{s_{\ell+2}}{\log q}<\frac{1}{|\Lambda|},
\end{equation}
where for notational convenience 
we put $\Lambda=\Lambda_{\ell+1}$.
We need to lower bound $|\Lambda|$.
Note that $ \vert p^{s_{\ell+1}} \cdot q^{-r_{\ell+1}}-1 \vert = \vert \exp(\Lambda)-1 \vert $. Either $ \vert \Lambda \vert \ge 1/2$, or $ \vert \Lambda \vert <1/2$, in which case 
$$
 \vert \exp(\Lambda)-1 \vert <2 \vert \Lambda \vert .
$$
We lower bound the left--hand side above using Theorem \ref{Matveev11}. We get
$$
2 \vert \Lambda \vert > \vert p^{s_{\ell+1}} \cdot p^{-r_{\ell+1}}-1 \vert > \exp\bigl(-1.4\cdot 30^5\cdot 2^{4.5} \,(1+\log \max\{s_{\ell+1},r_{\ell+1}\})\,\log p\cdot \log q\bigr).
$$
Note that we can assume that $s_{\ell+1}>r_{\ell+1}$. Indeed, otherwise we have $r_{\ell+1}\ge s_{\ell+1}$ and so
$$
 \vert \Lambda \vert = \vert (r_{\ell+1}-s_{\ell+1})\log q+s_{\ell+1}(\log q-\log p) \vert \ge \log\big(\frac{q}{p}\big)\ge \log\big(1+\frac{1}{p}\big)>\frac{1}{2p},
$$
which is a better inequality. 
In particular, $s_{\ell+1}\ge 2$. We thus get
\begin{align*}
 \vert \Lambda \vert &>\exp\bigl(-\log 2-1.4\cdot 30^5\cdot 2^{4.5} \,(1+\log s_{\ell+1}) \log p \cdot \log q\bigr)\\
&>\exp\bigl(-2\cdot 10^9\, \log s_{\ell+1} \,\log p\cdot \log q\bigr),
\end{align*}
where we have used the fact that $1+\log s_{\ell+1}\leq 2.5\log s_{\ell+1}$, for $s_{\ell+1}\geq 2$.
Thus, inequality \eqref{eq:2} gives
$$
\frac{s_{\ell+2}}{\log q}<\frac{1}{ \vert \Lambda \vert }<s_{\ell+1}^{C_1}.
$$
Hence, we get
$$
s_{\ell+1}>0.5\,s_{\ell+2}^{C_2}.
$$
Then, recalling 
\eqref{starting}-\eqref{cfrac-approx}, we obtain
\begin{equation}
\label{eq:Case2}
\log \Bigl(\frac{n_{i+1}}{n_i}\Bigr)
<
\frac{\log q}{s_{\ell+2}}
\le
\frac{\log q}{s_{\ell+1}}
<\frac{2\log q}{s_{\ell+2}^{C_2}}.
\end{equation}
We thus get that in all cases (\eqref{eq:Case1} and \eqref{eq:Case2}) we have
\begin{equation}
\label{eq:3}
\log\Bigl(\frac{n_{i+1}}{n_i}\Bigr)
<\frac{\log q}{s_{\ell+1}}<\frac{2\log q}{s_{\ell+2}^{C_2}}<\frac{2\log q}{u^{C_2}}.
\end{equation}
Since by assumption $p^{u}>q^v$ we have $p^{u}>n_i^{1/2}$ and hence
$u>\frac{\log n_i}{2\log p}$ we obtain
$$
\log\Bigl(\frac{n_{i+1}}{n_i}\Bigr)<\frac{2\,(2\log p)^{C_2}\,\log q}{(\log n_i)^{C_2}}<\frac{4\log q}{(\log n_i)^{C_2}}.
$$
But $(n_{i+1}-n_i)/n_i\in (0,p-1]$. In this interval, the image of the function 
$
\log(1+x)/x
$
is in the interval $[\frac{\log p}{p-1},1]$. Writing $n_{i+1}/n_i=1+x$, where $x=(n_{i+1}-n_i)/n_i$, we get that 
$$
n_{i+1}-n_i<\frac{4\,(p-1)\,n_i\,\log q}{(\log n_i)^{C_2}\,\log p},
$$
and we can choose $C_4$ to be at least 
$$
\frac{4(p-1)\,\log q}{\log p}.
$$
But, we also needed that this constant exceeds $q(\log q)^{C_2}$. Since $(\log q)^{C_2}<2$, it follows that 
if we choose $C_4\ge 4q$, then everything works out.

We proceed by a similar argument in the remaining case $q^v>p^{u}$. 
Note that $v\ge 1$. We take $\ell$ such that $r_{\ell}\le v< r_{\ell+1}$. We choose $\ell'\in \{\ell,\ell+1\}$
such that $r_{\ell'}/s_{\ell'}$ is smaller than $\theta_{p,q}$.
Note that if $\ell=0$, we then have $r_0/s_0=0<\theta_{p,q}$. So, $\ell'\ge 0$ is well-defined even if $\ell=0$. 
We consider the number
$$
p^{u+s_{\ell'}} \cdot q^{v-r_{\ell'}}, 
$$
which  is an integer which is a 
$\{p,q\}$-unit exceeding $n_i$, so it is at least $n_{i+1}$. Assuming $s_{\ell+1}\ge 2$
and reasoning as before we deduce \eqref{eq:3} 
and so 
\begin{equation}
\label{secondcase}
\log\Bigl(\frac{n_{i+1}}{n_i}\Bigr)<\frac{2\log q}{s_{\ell+2}^{C_2}}.
\end{equation}
The above inequality assumes that $s_{\ell+1}\ge 2$ which holds for all $\ell\ge 0$ except when $\ell=0$ and $q\in (p,p^2)$. But this is impossible since  we  must also have 
$0=r_0\le v<r_1=1$, so $v=0$, therefore $i=1$, which is false. 
We need to recast the above inequality in terms of $v$. By \eqref{cfrac-approx} we have
$$
 \vert r_{\ell+2}\log q-s_{\ell+2}\log p \vert <\frac{\log q}{s_{\ell+2}}.
$$
Hence, using that $s_{\ell+2} \ge s_{\ell+1}\ge 2$ we obtain
$$
s_{\ell+2}\log p\ge \Bigl(r_{\ell+2}-\frac{1}{s_{\ell+2}}\Bigr)\log q\ge 0.5\,r_{\ell+2}\log q.
$$
 Thus, we get
$$
s_{\ell+2}\ge  \frac{\log q}{2\log p}  r_{\ell+2},
\quad 
s_{\ell+2}^{C_2}
\ge  \Bigl(\frac{\log q}{2\log p}\Bigr)^{C_2}r_{\ell+2}^{C_2}>  \frac{r_{\ell+2}^{C_2}}2.
$$
By \eqref{secondcase} we get
$$
\log\Bigl(\frac{n_{i+1}}{n_i}\Bigr)<\frac{4\log q}{r_{\ell+2}^{C_2}}<\frac{4\log q}{v^{C_2}}.
$$
So, we get an upper bound in the right--hand side by a factor of at most $2$ larger than in the case when $p^{u}>q^v$. Following along we obtain that in this case
$$
n_{i+1}-n_i<C_4\frac{n_i}{(\log n_i)^{C_2}},
$$
we must have $C_4>8(p-1)
\frac{\log q}{\log p}$, and
so taking $C_4=8q$ suffices.

%&&&&&&&&&&&&&&&&&&&&&&&&&&&&&&&&&&&&&&&&&&&&&&%
\section{Bertrand's Postulate for the sequence \texorpdfstring{$(n_i)_{i\ge 0}$}{ni}} \label{sec3}
%&&&&&&&&&&&&&&&&&&&&&&&&&$$$$$$$$$$$$$$$$$$$$$$$$$$%

%---------------------------------------------------------------------
\subsection{The statement}
\label{sec:statement}
%---------------------------------------------------------------------
We will show a way to compute $n_{p,q}(\alpha)$. 
Tables \ref{Table-n2q-5/3} and
\ref{Table-n2q-3/2} give some examples for $p=2$ and $\alpha=5/3$, respectively $3/2$, 
some of which are relevant for an application.
\par Basic results from the
theory of Diophantine 
approximation
(cf.\,Section \ref{sect:comp_npq_alpha}), 
ensure the existence of integers $e,f,g$, and $h$ such that
$$1<\frac{p^f}{q^e}<\alpha\quad \textrm{and} \quad 1<\frac{q^h}
{p^g}<\alpha.$$
We claim that $n_{p,q}(\alpha)\le  p^g\cdot q^e.$ In order to see this, one can observe
that any integer $n:= p^\ell \cdot q^k\ge p^g\cdot q^e$ satisfies either
$k\ge e$, or $\ell\ge g.$ In case $k\ge e,$ we note that the
number $ p^{\ell+f}\cdot q^{k-e}$ lies in $[n,n\alpha).$
In case $\ell\ge g,$ we have $p^{\ell-g}\cdot q^{k+h}\in [n,n\alpha).$ Next 
one tries to find
an integer
$n_{\text{new}}:= p^\ell\cdot q^k\in [\lfloor  p^g\cdot q^e/\alpha \rfloor, p^g\cdot q^e)$, where
$\lfloor x\rfloor$ denotes the integral part of $x$. 
If successful, we continue until we 
fail, each time considering the interval
$[\lfloor n_{\text{new}}/\alpha\rfloor,n_{\text{new}})$.
In Section \ref{sect:comp_npq_alpha} we present a much more refined way of determining $n_{p,q}(\alpha)$.

\subsection{An upper bound for \texorpdfstring{$n_{p,q}(\alpha)$}{npqalpha}}
Recall that the \emph{discrepancy} $D_N$ of a sequence $(a_m)_{m=1}^N$
of real numbers (not necessarily distinct) is defined as
$$
D_N=\sup_{0\le \gamma\le 1}\Bigr|\frac{\#\{m\le N~:~\{a_m\}<\gamma\}}{N}-
\gamma\Bigr|,
$$
where $\{x\}$ denotes the fractional part of a real number $x$.
The Koksma-Erd\H os-Tur\'an inequality (see, for 
example, Kuipers and Niederreiter \cite[Lemma 3.2]{KuiNied}) states that
\begin{equation}
\label{eq:ErdosTuran}
D_N\le \frac{3}{H} + \frac{3}{N}\sum_{m=1}^{H}\frac{1}{m\|a_m\|},
\end{equation}
where $\|x\|$ is the distance from $x$ to the nearest
integer and $H\le N$ is an arbitrary positive integer.
In this section we will improve on the bound \eqref{directbound} by applying 
this inequality
to upper bound the discrepancy of the sequence $(j\theta_{p,q})_{j\ge 1}$, with 
$\theta_{p,q}:=\log p/\log q$. 
\begin{theorem}
\label{lem:10}
Let $p<q$ be primes and $\alpha\in (1,p]$.
Put $C_5=10^9\log p\cdot \log q$.
There are positive integers $f$ and $g$ such that 
\begin{equation}
\label{eq:containment}
\Bigl\{f\frac{\log q}{\log p}\Bigr\}\in \Bigl(0,\frac{\log\alpha}{\log q}\Bigr), 
\qquad 
\Bigl\{ g\frac{\log p}{\log q}\Bigr\}\in \Bigl(1-\frac{\log \alpha}{\log q},1\Bigr),
\end{equation}
and 
\begin{equation}
    \label{f+gbound}
\max\{f,g\}<\Bigl(\frac{60\log q}{\log \alpha}\Bigr)^{C_5}.
\end{equation}
In particular,
$$
n_{p,q}(\alpha)\le q^e\cdot  p^g<p^{f+g}
<
\exp\Bigl(2\log p\Bigl(\frac{60\log q}{\log \alpha}\Bigr)^{C_5}\Bigr).
$$
\end{theorem}
\begin{Proof}
From the definition of 
$D_N$ we see that the inequality
\begin{equation}
\label{ineqbasic}
\#\{m\le N~:~\alpha\le \{a_m\}<\beta\}\ge (\beta-\alpha)N-2D_N N
\end{equation}
holds for all $0\le \alpha\le \beta\le 1$.
We will apply this with $a_m=m\,\theta_{p,q}$ for all $m=1,\ldots,N$. Writing
$$
I=\Bigl(0,\frac{\log\alpha}{\log q}\Bigr),\qquad J=\Bigl(1-\frac{\log \alpha}{\log q},1\Bigr),
$$
both intervals of length $\log \alpha/\log q$, it follows from \eqref{ineqbasic} that
\begin{equation}
\label{eq:ineqforN}
\#\{m\le N: \{a_m\}\in I\}  \ge   | I |  N
-2D_N N=   \Bigl(\frac{\log \alpha}{\log q}\Bigr)N-2D_N N,
\end{equation}
and similarly with $I$ replaced by $J$. 
In particular, if the right-hand side is positive, then there is $u\le
N$ with $\{a_u\}\in I$.
We will now upper bound $D_N$
using \eqref{eq:ErdosTuran}.
To bound $\|a_m\|$, note that
$$
\|a_m\|=\Bigl|m\frac{\log p}{\log q}-t\Bigr|= \frac{|\Lambda|}{\log q},\text{~with~}
\Lambda:=m\log p-t\log q,
$$
for some integer $t$ with $t\le m(\log p)/(\log q)+1<2m$.  
Since $\theta_{p,q}$ is irrational we
have either $|\Lambda|\ge 1/2$, or $0<|\Lambda|<1/2$. If $0<|\Lambda|<1/2$, then
$$
2|\Lambda|>|\exp(\Lambda)-1|,
$$
and we can apply Matveev's theorem to get a lower bound 
on $|\exp(\Lambda)-1|$ and hence on $|\Lambda|$.
In both cases we take in Matveev's theorem
$$
k=2,\quad \gamma_1=p,\quad \gamma_2=q,\quad b_1=m,\quad b_2=-t,
$$
and, noting that we can set $B:=2m$, we get
\begin{equation}
\label{eq:Matv}
2|\Lambda|  > \exp\left(-c_1(\log p) (1+\log (2m))  \log q\right),
\end{equation}
where $c_1=1.4\cdot 30^5\cdot 2^{4.5}<8\cdot 10^8-\log 2$. We get
\begin{eqnarray*}
|\Lambda| & > & \exp(-8\cdot 10^8(1+\log(2m))
\log p\cdot \log q)\\
& = &
q^{-8\cdot 10^8
(1+\log (2m))\log p}\quad {\text{\rm for}}\quad m\ge  1.
\end{eqnarray*}
We thus obtain that, if $H\ge 60$ and $2m\le H,$ then
$$
1+\log (2m)\le 1+\log H\le 1.245\log H\qquad (H\ge 60),
$$
and so inequality \eqref{eq:Matv} leads to
\begin{eqnarray*}
\frac{1}{\|a_m\|} & \le & (\log q)\,q^{(8\cdot 1.245) \cdot 10^8 \log
H\log p}<q^{(10^9-3)\log H\log p}\\
& = & H^{(10^9-3)\log p\cdot \log q}<
H^{C_5-2}.
\end{eqnarray*}
Thus,
$$
D_N\le 
3\Bigl(\frac{1}{H}+\frac{H^{C_5-2}}{N} \sum_{m=1}^H \frac{1}{m} \Bigr)
<3\Bigl(\frac{1}{H}
+\frac{H^{C_5-1}}{N}\Bigr),
$$
where we trivially bounded the sum by $H$.
Choosing $H:=\left\lfloor N^{1/C_5}\right\rfloor$ we get,
assuming
still that $H\ge 60$ and therefore that
\begin{equation}
\label{eq:hyp}
N^{1/C_5}\ge 60,\quad {\text{\rm which is equivalent
to}}\quad N\ge 60^{C_5},
\end{equation}
we obtain
$$
D_N \le 3\Big(\frac{1}{H}+\frac{H^{C_5-1}}{N}\Big) \le  
3\Bigl(\Big\lfloor N^{1/C_5}\Big\rfloor^{-1}+N^{-1/C_5}\Bigr)\\
 \le 7N^{-1/C_5}.
$$
To derive the final inequality we used the
trivial observation that if $x\ge 60$, then
$$
\frac{1}{\lfloor x\rfloor}+\frac{1}{x}\le 
\frac{1}{x}\Big(\frac{1}{1-\frac{1}{x}}+1\Big)
\le\frac{119}{59}\cdot\frac{1}{x}<\frac{7}{3}\cdot 
\frac{1}{x}.$$
Turning now our attention to the inequality \eqref{eq:ineqforN}, we see that
\begin{equation}
\label{koksma}    
N\Bigl(\frac{\log \alpha}{\log q}-2D_N\Bigr)>
N\Bigl(\frac{\log \alpha} {\log q}
-14N^{-1/C_5}\Bigr).
\end{equation}
Thus, if $N\ge N_0$ with
\begin{equation}
\label{eq:N0}
N_0:= \Bigl(\frac{60\log q}{\log \alpha}\Bigr)^{C_5},
\end{equation}
and hence $N_0\ge 60^{C_5}$, the right--hand side of 
\eqref{koksma} is at 
least
\begin{equation}
\label{N15}
\frac{23}{30}\,N\,\frac{\log \alpha}{\log q} 
\end{equation}
and hence positive.
From what we have seen at the beginning of Section \ref{sec:statement},
we have  $n_{p,q}(\alpha)\le q^e\cdot  p^g<p^{f+g}$. 
The proof is now completed on invoking 
\eqref{f+gbound}.
\end{Proof}
A similar proof also appears in
Ferrari et al.\,\cite{FerLM} in the
context of bounding the largest exceptional value of the discriminator of certain Lucas sequences.
\subsection{Application to the discriminator of Lucas type sequences}
\label{sec:lucasdisc}
Faye et al.\,\cite{FLM} considered the sequence $\{U_n(k)\}_{n\ge 0}$ defined uniquely by
$$
U_{n+2}(k)=(4k+2)U_{n+1}(k)-U_n(k),~~U_0(k)=0,~U_1(k)=1.
$$
For $k = 1$, the sequence
is 
$$
0, 1, 6, 35, 204, 1189, 6930, 40391, 235416, 1372105,\ldots,
$$
which is A001109 in OEIS. On noting that
$$U_{n+2}(k) - U_{n+1}(k)
= 4kU_{n+1}(k) + U_{n+1}(k) -
U_n(k)\ge 1,$$
one sees that the sequence
$U_n(k)$ consists of strictly increasing non-negative numbers.
We can now define its \emph{discriminator} ${\mathcal D}_k(n)$ as the 
smallest positive integer $m$ such that $U_{0}(k),\ldots,U_{n-1}(k)$ are pairwise 
distinct modulo $m$. 
\par Primes of the form $2^n-1$ are
called \emph{Mersenne primes}. Note that $n$ has to be a prime.
The first few Mersenne primes are: 
$3, 7, 31, 127, 8191, 131071, 524287,\ldots$
The first few with
exponent $n\equiv 1\bmod{4}$ are 
$31,8191, 131071, 2305843009213693951,\ldots$ Choosing
$k$ to be a Mersenne prime, it turns out that
${\mathcal D}_k(n)$ is a $\{2,k\}$-unit for all $n$ large enough.
\begin{proposition}
\label{prop:disc}
Let $p$ be an odd prime such that also $q:=2^p-1$
is a prime number. Let $(n_i)_{i\geq 0}$ be the sequence of
consecutive $S$-units with $S=\{2,q\}$.
Then 
\begin{equation}
    \label{discineq}
{\mathcal D}_q(n)\le \min\{n_i\ge n\},
\end{equation}
with equality if the interval $[n,3n/2)$ contains
an $S$-unit $n_j$. In particular, if
$n\ge n_{2,q}(3/2)$ we have equality in \eqref{discineq}.
If $p\equiv 1\bmod{4}$ and $p>5$, then we have 
equality for $n\ge n_{2,q}(5/3)$ and
${\mathcal D}_q(n)$ is a $\{2,5,q\}$-unit for every $n\ge 1$.
\end{proposition}

\begin{Proof}
The first assertion follows on taking for $k$ the Mersenne prime $q$ in Theorem 3 
of Faye et al.\,\cite{FLM}, the second 
assertion is a consequence of 
Theorem 3 of Ferrari et al.\,\cite{FLM}.
\end{Proof}

\renewcommand{\arraystretch}{1.2}
\begin{table}[htp] 
\scalebox{0.9}{
\begin{tabular}{ | r | r |}
\hline
$q$\phantom{0} & $n_{2,q}(5/3)$ \\ \hline
$3$ & $1$ \\ \hline 
$5$ & $1$ \\ \hline 
$7$ & $1$ \\ \hline 
$11$ & $1$ \\ \hline 
$13$ & $1$ \\ \hline 
$17$ & $39322$ \\ \hline 
$19$ & $154$ \\ \hline 
$23$ & $10$ \\ \hline 
\end{tabular}
\hskip0.25cm
\begin{tabular}{ | r | r |}
\hline
$q$\phantom{0} & $n_{2,q}(5/3)$ \\ \hline
$29$ & $279$ \\ \hline 
$31$ & $274839850$ \\ \hline 
$37$ & $615$ \\ \hline 
$41$ & $20$ \\ \hline 
$43$ & $20$ \\ \hline 
$47$ & $20$ \\ \hline 
$53$ & $20$ \\ \hline 
$59$ & $66836$ \\ \hline 
\end{tabular}
\hskip0.25cm
\begin{tabular}{ | r | r |}
\hline
$q$\phantom{0} & $n_{2,q}(5/3)$ \\ \hline
$61$ & $4358036$ \\ \hline 
$67$ & $10066330$ \\ \hline 
$71$ & $2458$ \\ \hline 
$73$ & $2458$ \\ \hline 
$79$ & $39$ \\ \hline 
$83$ & $39$ \\ \hline 
$89$ & $39$  \\ \hline
$97$ & $39$  \\ \hline
\end{tabular}
} 
\caption{\small{$n_{2,q}(5/3)$; $3\le q\le 100$.}}
\label{Table-n2q-5/3}
\end{table} 

\renewcommand{\arraystretch}{1.2}
\begin{table}[htp] 
\scalebox{0.9}{
\begin{tabular}{ | r | r |}
\hline
$q$\phantom{0} & $n_{2,q}(3/2)$ \\ \hline
$3$ & $1$ \\ \hline 
$5$ & $11$ \\ \hline 
$7$ & $131$ \\ \hline 
$11$ & $6$ \\ \hline 
$13$ & $70$ \\ \hline 
$17$ & $699051$ \\ \hline 
$19$ & $171$ \\ \hline 
$23$ & $11$ \\ \hline 
\end{tabular}
\hskip0.25cm
\begin{tabular}{ | r | r |}
\hline
$q$\phantom{0} & $n_{2,q}(3/2)$ \\ \hline
$29$ & $8971$ \\ \hline 
$31$ & $282022636380491$ \\ \hline 
$37$ & $683$ \\ \hline 
$41$ & $683$ \\ \hline 
$43$ & $22$ \\ \hline 
$47$ & $22$ \\ \hline 
$53$ & $1131$ \\ \hline 
$59$ & $4381419$ \\ \hline 
\end{tabular}
\hskip0.25cm
\begin{tabular}{ | r | r |}
\hline
$q$\phantom{0} & $n_{2,q}(3/2)$ \\ \hline
$61$ & $18018054422$ \\ \hline 
$67$ & $2932031007403$ \\ \hline 
$71$ & $174763$ \\ \hline 
$73$ & $174763$ \\ \hline 
$79$ & $2731$ \\ \hline 
$83$ & $2731$ \\ \hline 
$89$ & $43$  \\ \hline
$97$ & $4139$  \\ \hline
\end{tabular}
} 
\caption{\small{$n_{2,q}(3/2)$; $3\le q\le 100$.}}
\label{Table-n2q-3/2}
\end{table} 
\par Numerical work by Matteo Ferrari
never led to inequality in \eqref{discineq}
in case $q=2^p-1$ with $p\equiv 3\bmod{4}$, but in 
case $p\equiv 1\bmod{4}$ he found examples, 
e.g.,
$${\mathcal D}_{8191}(129)=250, 
\quad {\mathcal D}_{131071}(129)=250, \quad {\mathcal D}_{131071}(257)=500.$$ As predicted by Proposition \ref{prop:disc}, these values are indeed $\{2,5,2^p-1\}$-units. 
\par One finds that $n_{2,7}(3/2)=131$ and computing ${\mathcal D}_7(n)$ for the integers $1\le n\le 130$ shows that 
for $q=7$ there is always equality in \eqref{discineq}, in line with the observations of Ferrari. 
Unfortunately, we are not able
to classify ${\mathcal D}_q(n)$ completely for any Mersenne prime $q>7$, as
then $n_{2,q}(5/3)$ appears to be very
large (and hence $n_{2,q}(3/2)$ even
more so). The Mersenne prime $127$, 
for example, leads to a very large value 
of $n_{2,127}(3/2)$:
the algorithm, which we will 
describe in Section \ref{sect:comp_npq_alpha},
starts with the pair of exponents $(615, 1)$ 
and, after having generated $23\,150$ left neighbors, 
terminates with the pair $(6,36)$.
Hence, the final result is 
$n_{2,127}(3/2)   
= \lceil 2^6 \cdot 127^{36}\cdot 2/3 \rceil$,
a number having $78$ digits.
\par For $n_{2,8191}(5/3)$ the algorithm starts with the pair
$(73\,800, 1)$ and, after having generated $195$ $078\,401$ left neighbors, 
terminates with the pair $(12,1\,493)$.
Hence, the final result is 
$n_{2,8191}(5/3) 
=
\lceil 2^{12} \cdot 8191^{1493}\cdot 3/5 \rceil
$,
a number having $5\,847$ digits.
The whole computation of this case took about $20$ minutes and $30$ seconds.

\par For $n_{2,131071}(5/3)$ the algorithm starts with the pair
$(1\,544\,466, 1)$, so that $n_{\textrm{start}}$
has $464\,936$ digits.
After having generated $65\,305\,850\,756$ left neighbors
(more than $65$ billion neighbors \dots), 
it terminates with the pair $(16,23\,897)$. 
Hence, the final result is 
$n_{2,131071}(5/3) 
=
\lceil 2^{16} \cdot 131071^{23897}\cdot 3/5 \rceil
$,
a number having $122\,298$ digits!!
The whole computation of this case took about $5$ days and $17$ hours.
The running times were obtained
on a Dell Optiplex $3050$, Intel i$5$-$7500$, $3.40$GHz,
$32$GB of RAM running Ubuntu $22.04.03$-LTS and PARI/GP $2.15.4$.

Apart from having  sharper estimates for $n_{\textrm{start}}$, the
 only hope to obtain $n_{2,q}(3/2)$ or 
$n_{2,q}(5/3)$ in case $q > 131071$ is a Mersenne prime,
is to implement part of this algorithm using the C programming language. This is doable
since, once one has gained the knowledge of a sufficiently large number of the
$\theta_{2,q}$-convergents, the continued fraction part and the first trivial part
can work on exponents only. On the other hand, for the second trivial part the 
use of multiprecision arithmetic is mandatory.

\renewcommand{\arraystretch}{1.2}
\begin{table}[htp] 
\scalebox{0.9}{
\begin{tabular}{ | r | r | r | r | r |}
\hline
$p$ & $q=2^p-1$ & $\alpha$ & $n_{2,q}(\alpha)$\phantom{012345} & $\ell(n_{2,q}(\alpha))$ \\ \hline
$2$ & $3$ & $3/2$ & $1$ & $1$\\ \hline 
$3$ & $7$ & $3/2$ & $131$  & $3$\\ \hline 
$5$ & $31$ & $5/3$ & $274\,839\,850$  & $9$\\ \hline 
$7$ & $127$ & $3/2$ &   $ \lceil 2^6 \cdot 127^{36}\cdot 2/3 \rceil$ & $78$\\ \hline 
$13$ & $8\,191$ & $5/3$ &  $\lceil 2^{12} \cdot 8191^{1493}\cdot 3/5 \rceil$ & $5\,847$\\ \hline 
$17$ & $131\,071$ & $5/3$ & $\lceil 2^{16} \cdot 131071^{23897}\cdot 3/5 \rceil$ & $122\,298$ \\ \hline 
\end{tabular}
} 
\caption{Known values of $n_{2,M_p}(\alpha)$ for Proposition \ref{prop:disc}; $\ell(n)$ 
is the number of decimals of $n$.}
\label{Table-Mersenne}
\end{table} 

%&&&&&&&&&&&&&&&&&&&&&&&&&%
 \section{Computations}\label{sec4}
%&&&&&&&&&&&&&&&&&&&&&&&&&%
In this section, we will use heavy computations
 to search for $n_{p,q}(\alpha)$ and to study
the extremal and average behaviors of the gaps $n_{i+1}-n_i$.
In both cases, we will use the knowledge of the continued fractions
convergents of $\theta_{p,q}$ to generate left or right neighbors
of a given $S$-unit $n_i$.

The first computational problem we address is to search for $n_{p,q}(\alpha)$
since, as we will see later, its knowledge plays also a role in studying
the extremal and average behaviors of the gaps $n_{i+1}-n_i$.

\subsection{Computation of \texorpdfstring{$n_{p,q}(\alpha)$}{npqalpha}}
\label{sect:comp_npq_alpha}
We now explain how the algorithm to compute $n_{p,q}(\alpha)$ works.
First of all, thanks to  Theorem \ref{lem:10}, we know that it is possible to 
identify $n_{\textrm{start}}$, a suitable an upper bound for $n_{p,q}(\alpha)$.
After having identified such a point, we will 
then generate a part of the $n_i$ sequence by looking for $S$-units less than $n_i$.

We will now describe a way how to generate such a suitable starting
point.

\subsubsection{Generating the starting point \texorpdfstring{$n_{\textrm{start}}$}{nstart}}
\label{nstart-gen}
An easy argument, see 
Section \ref{sec:statement}, shows that if $f,e,g,h\ge 0$ are such that 
\begin{equation}
\label{twoinequalities}
1<\frac{p^f}{q^e}<\alpha\quad \textrm{and} \quad 1<\frac{q^h}{p^g}<\alpha
\end{equation}
then $n_{p,q}(\alpha)\le p^g \cdot q^e$. 
On taking logarithms these inequalities can be rewritten as
$$
0<\theta_{p,q}-\frac{e}{f}<\frac{\log \alpha}{f\log q}
\qquad \textrm{and} \qquad 
0<\frac{h}{g}-\theta_{p,q}<\frac{\log \alpha}{g\log q}.
$$
These type of inequalities appear in the theory of continued fractions and suggest to take
\begin{equation}
\label{doublefrac}    
\frac{e}f=\frac{r_k}{s_k}
\quad\textrm{and}\quad\frac{h}{g}=
\frac{r_{\ell}}{s_{\ell}},
\end{equation}
where $k$ is even and minimal, $\ell$ is odd and minimal and $(r_m/s_m)_{m\ge 0}$ is the sequence of convergents of $\theta_{p,q}$. Using the inequality
$$
\Bigl\vert \theta_{p,q}-\frac{r_m}{s_m} \Bigr\vert <\frac{1}{s_m s_{m+1}},
$$
see, e.g., \cite[Thm.\,171]{HW},
we deduce that in both cases  ($m=k$ and $m=\ell$) it suffices to require that 
$$
\frac{1}{s_m s_{m+1}}<\frac{\log \alpha}{s_m \log q},
$$
that is
\begin{equation}
\label{m-condition}
s_{m+1}>\frac{\log q}{\log \alpha}.
\end{equation}
Clearly, the larger $q$ 
will be, or the closer to $1$ the value of $\alpha$ will be, 
the larger $m$ will become.
Let $M$ be the minimal $m\ge 0$ such that \eqref{m-condition} holds.
In practice, to determine $M$ we need a way to generate the continued fractions convergents
of $\theta_{p,q}$. To do this, we heavily rely on the PARI/GP \cite{PARI2023} internal functions.
Letting $\eps\in(0,1)$ be the required accuracy for the computations, the \texttt{contfrac}$(x)$
function of PARI/GP returns the list of the partial quotients $[a_0,\dotsc,a_n]$ 
of the continued fraction expansion of $x>0$ so that $\vert x - (a_0+1/(a_1+\dotsm+1/a_n)) \vert  < \eps$.
Using such a continued fraction expansion of $x$, the \texttt{contfracpnqn} function of PARI/GP returns 
two sorted lists containing $p_n\in {\mathbb N}^*$ and $q_n\in {\mathbb N}^*$, $(p_n,q_n)=1$, $n\in {\mathbb N}$,
the numerators and the denominators of the convergents $p_n/q_n$ of $x$.
In this way, it is relatively easy to obtain both the upper and lower convergents
for each fraction $\theta_{p,q}$, $2\le p<q$ both primes, we have to work with.
Our practical computations are performed with $\eps=10^{-19}$.

Having now such a sorted list of the denominators $s_{m}$ of the convergents for $\theta_{p,q}$,
 a standard dyadic search procedure will quickly provide
$M$.
Nevertheless, it is possible to obtain some theoretical information
about $M$; this might be useful in the case one has initially no access to the sorted list 
of the convergent denominators.
Since $s_{m+1}\ge F_{m+2}$, where $F_k$ is the $k$-th member of the Fibonacci sequence, it is enough to require that $F_{m+2}\ge \log q/\log \alpha$.  
Recalling that $F_{m+2}\ge \phi^m$, 
where $\phi=(1+{\sqrt{5}})/2$, it suffices to take
\begin{equation}
\label{startestimate}    
m\ge \frac{\log\log q-\log\log \alpha}{\log \phi} >0.
\end{equation}
Hence $M \le \lceil \frac{\log\log q-\log\log \alpha}{\log \phi} \rceil$.
At the cost of complicating our inequality we can actually do a bit better.
This improvement often turns out to make a real difference in numerical practice and so is worth the effort; this is due to the fact that both the numerators 
and the denominators of the continued fractions convergents are 
exponentially fast increasing 
sequences. As a consequence, being able to choose a smaller value for $M$
ensures that much
smaller values for the exponents in the definition
of $n_{\textrm{start}}$ can be chosen, see \eqref{nstart-def}.

Namely, we will use the sharper inequality 
$F_{k}> \frac{\phi^k}{\sqrt{5}} - \frac{53}{500}$ for every $k\ge 3$
(this inequality is a consequence of Binet's formula $F_{k}=\frac{\phi^k-(-\phi)^{-k}}{\sqrt{5}}$
and the observation that for
$k\ge 3$ one has $\vert\frac{(-\phi)^{-k}}{\sqrt{5}}\vert < \frac{53}{500}$).
Hence if 
$$
m + 2  > 
\frac{\log\bigl(\frac{\log q}{\log \alpha} + \frac{53}{500}\bigr)+\frac{\log 5}{2}}{\log \phi} 
$$
for every $m \ge 1$, then $F_{m+2}\ge \log q/\log \alpha$.
It is not hard to verify that
$\frac{\phi^{m+2}}{\sqrt{5}} - \frac{53}{500} > \phi^m$
for every $m \ge 0$,
%for every $m \ge 1$, 
and so the latter inequality gives a better lower bound
for $m$ than \eqref{startestimate}, namely
$$
m \ge
\Bigg\lceil  \frac{\log\bigl(\frac{\log q}{\log \alpha} 
+ \frac{53}{500}\bigr)+\frac{\log 5}{2}}{\log \phi} \Bigg\rceil  - 2
\ge M.
$$

Once $M$ is determined with the dyadic search procedure or
using the previously described estimates, 
we can take
%So, taking 
$\{k,\ell\}=\{M,M+1\}$ in \eqref{doublefrac},
and get that
\begin{equation}
\label{nstart-def}
n_{p,q}(\alpha)\le p^{g} \cdot q^{e} \le n_{\text{start}}: = \min(p^{s_{M+1}}\cdot q^{r_M}; p^{s_{M}}\cdot q^{r_{M+1}}).
\end{equation}
\par In \eqref{twoinequalities} we want to find the solutions with
$p^g$ and $q^e$ minimal. It follows from the
basics of continued fractions that the approach with the convergents as described here is actually optimal. So in retrospect our choice in \eqref{doublefrac} was best possible.

\subsubsection{Generating left neighbors}
\label{left-neighbors-gen}
We take $n_{\text{start}}$ as starting 
candidate and proceed as follows.   
Assume that we know $n_i = n_{\text{start}}$ 
and we want to obtain an $S$-unit $n^{*}_{i-1}\le n_{i-1}$  such that
$n^{*}_{i-1}  \ge \lfloor n_i/\alpha \rfloor$. Remark that the  goal 
here is not to find $n_{i-1}$, the predecessor element in the $n_i$ sequence, but just an
$S$-unit less than $n_i$ that verifies the Bertrand Postulate condition. 
For this reason in this procedure we will always choose, if possible, the 
$S$-unit $n^{*}_{i-1}$ having the maximal distance from $n_i$ 
compatible with the condition $n^{*}_{i-1} \ge \lfloor n_i/\alpha \rfloor$, because this reduces the total amount of
computations to be performed to determine $n_{p,q}(\alpha)$.

In order to achieve this, we combine three different ways of searching for 
$n^{*}_{i-1}$ in the following algorithm.
\begin{enumerate}[label={\alph*)}, wide, labelindent=0pt, nosep, after=\vspace{10pt}]
\item 
\label{cfrac-search}
\textbf{Searching using continued fraction convergents.}
Letting  $n_i:= p^a \cdot q^b$, we have two possible choices.
\begin{enumerate}[label={\arabic*)}, leftmargin=*]
\item We search for $\ell$ such that $a\in [s_{\ell}, s_{\ell+1}]$,
where $(r_m/s_m)_{m\ge 0}$ is the sequence of convergents to $\theta_{p,q}$. 
Choose $l$ in $\{\ell,\ell-1\}$ such that $r_l/s_l<\theta_{p,q}$;
this is equivalent to $x_1:= - s_l \log p + r_l \log q <0$
and hence $q^{r_l}/p^{s_l} < 1$.
Define  $N_1:=n_i\exp(x_1)< n_i$. 

\item  
We search for $\ell$ such that $b \in [r_{\ell},r_{\ell+1}]$.
Choose $l$ in $\{\ell,\ell-1\}$ such that
$r_l/s_l>\theta_{p,q}$;
this is equivalent to $x_2 := s_l \log p - r_l \log q < 0$
and hence $p^{s_l}/q^{r_l} < 1$.
Define $N_2:=n_i\exp(x_2)< n_i$.
\end{enumerate}
Both $N_1$ and $N_2$ are less than $n_i$ and thus candidates to be chosen as $n^{*}_{i-1}$.
We need now select the best of them;  i.e., the one whose distance from $n_i$ is maximal,
compatible with the condition $n^{*}_{i-1} \ge n_i/\alpha \ge \lfloor n_i/\alpha \rfloor$.
We point out that we compare here with $n_i/\alpha$, rather than 
with $\lfloor n_i/\alpha \rfloor$, 
since $n_i/\alpha$ allows us, by taking logarithms, to work on the exponents only; in other words,
we try to avoid as long as possible the necessity of performing
the costly computations of $n_i$ and $\lfloor n_i/\alpha \rfloor$.

Let now $x_{\textrm{max}} := \max(x_1;x_2)$ and $x_{\textrm{min}} := \min(x_1;x_2)$.
If $x_{\textrm{max}} < - \log \alpha$, then  $N_1,N_2$ are both $<n_i/\alpha$.
In this case we terminate this step with $n_i$ and we proceed with
step \ref{trivial-first-search};
we also remark that at this point we know that $n_{p,q}(\alpha)  \le  \lceil n_{i}/\alpha \rceil$.

Assume that $x_{\textrm{max}} \ge - \log \alpha$; this means that at least one of $x_1$ and $x_2$ is $\ge - \log \alpha$. If both $x_1$ and $x_2$ are $\ge - \log \alpha$, 
the best choice is to select the smallest one, $x_{\textrm{min}}$,
since $n_i\exp(x_{\textrm{min}})$ has the largest distance from $n_{i}$.
Hence,
if $x_{\textrm{min}}\ge - \log \alpha$, we choose $n^{*}_{i-1} = n_i\exp({x_{\textrm{min}}}) = \min(N_1;N_2)$; 
otherwise we are in the case $x_{\textrm{min}}<- \log \alpha \le x_{\textrm{max}}$, 
and we are forced to choose $n^{*}_{i-1} = n_i\exp({x_{\textrm{max}}})= \max(N_1;N_2)$.
In both cases, we replace $n_i$ with $n^{*}_{i-1}$  and we repeat step \ref{cfrac-search}.
 
\item
\textbf{Searching trivially: first part.}
\label{trivial-first-search}
Assuming that step \ref{cfrac-search}
terminates with $n_{i}$, an improved approximation of 
$n_{p,q}(\alpha)$ is then obtained by trivially searching for a value of
$x := a \log p + b \log q$
in the range $[\log (n_{i}/\alpha), \log n_{i})$,
where $a,b\in {\mathbb Z}$ run in the intervals
$0\le b \le \lfloor\log n_{i}/\log q\rfloor$ and
$0\le a \le \lfloor(\log n_{i} - b \log q)/\log p\rfloor$.
Moreover, we can also exploit the fact that the
search procedure of step \ref{cfrac-search} produced
$\overline{n} = \max(N_1;N_2) = p^{\overline{u}}\cdot q^{\overline{v}}$
such that $\overline{n} <  n_i/\alpha$. This means
that the value 
$\overline{x} := \overline{u} \log p + \overline{v} \log q$ is smaller than 
$\log( n_i/\alpha).$
Since we need to find $x \ge \log( n_i/\alpha)> \overline{x}$,
we must have either $a > \overline{u}$ 
or $b > \overline{v}$.
As a consequence, since $\log q> \log p$, the best strategy to determine $x$ 
is to first search for $\overline{v} + 1 \le b \le \lfloor\log n_{i}/\log q\rfloor$ and
$0\le a \le \lfloor(\log n_{i} - b \log q)/\log p\rfloor$.
If we have no success, we then work with
$0 \le b \le \overline{v}$  and
$\overline{u}+1 \le a \le \lfloor(\log n_{i} - b \log q)/\log p\rfloor$.

If in one of the previously described procedures we find 
a solution $x \ge\log( n_i/\alpha)$, we have obtained an 
$S$-unit $n^{*}_{i-1}:=\exp(x)$
such that $n_i/\alpha \le 
n^{*}_{i-1} < n_i$. 
In this case,
we replace $n_i$ with $n^{*}_{i-1}$ and we 
start again the search described in step \ref{cfrac-search}. 
If we do not find any solution, we continue with step \ref{trivial-second-search}.

\item
\textbf{Searching trivially: second part.}
\label{trivial-second-search}
After steps \ref{cfrac-search}-\ref{trivial-first-search} are over,
we have not yet determined 
$n_{p,q}(\alpha)$, since before, for efficiency reasons, 
we have replaced the condition
$\log ( \lfloor n_{i}/\alpha \rfloor ) \le x$
with the sharper, but easier to compute, 
$\log (n_{i}/\alpha) \le x$.
Hence we perform here another trivial search like 
the previous one, but using the the correct lower
bound for $x$ mentioned before.
If we find such a solution $x$, we have obtained an 
$S$-unit $n^{*}_{i-1}:=\exp(x)$
such that $ \lfloor n_{i}/\alpha \rfloor \le n^{*}_{i-1} < n_i$. In this case
we replace $n_i$ with $n^{*}_{i-1}$ and we 
start again the search described in step \ref{cfrac-search}.
If we do not find any solution $x$, this step has determined 
$n_0$, the smallest generated $S$-unit such that $n_{i-1} \ge \lfloor n_i/\alpha \rfloor$
holds for every $i\ge 1$.
In this case, we continue with step \ref{final-step}.

\item 
\textbf{Final computation.}
\label{final-step}
In the final step we obtain $n_{p,q}(\alpha) = \lceil n_0/\alpha \rceil$.
\end{enumerate}

The search in step \ref{cfrac-search} is clearly the fastest
one. Hence, the previously described procedure  optimizes
the computational cost by minimizing the number of times we are using the much slower
trivial searches of steps \ref{trivial-first-search}-\ref{trivial-second-search}.
Moreover, except
for step \ref{trivial-second-search}, in which case the presence of the floor function forces us to compute $n_i$,
the computations can be directly performed on 
the exponents $a,b$, rather than with the
prime powers involved. 
This requires far less memory usage and, at the same time,
as much smaller numbers are involved,
improves the running time of the algorithm.
We ran our program for $\alpha \in \{2, 7/4, 5/3, 3/2, 4/3\}$, 
$2\le p < q \le 500$, 
on the cluster\footnote{We used a machine having
$2$ x Intel(R) Xeon(R) CPU E$5$-$2630$L v$3$\@$1.80$GHz
and $192$GB of RAM (but we used up to $128$GB of RAM in our computations).} located at the Dipartimento di Matematica ``Tullio Levi-Civita'' of 
the University of Padova;
the running times were respectively
$6$ hours and $52$ minutes ($\alpha =2$),
$9$ hours and $39$ minutes ($\alpha =7/4$),
$11$ hours and $23$ minutes  ($\alpha =5/3$),
one day, $14$ hours and $36$ minutes ($\alpha =3/2$),
and $14$ days, $20$ hours ($\alpha = 4/3$).

To show the importance of working on the exponents only, 
we report here some data about  the computations for 
$n_{3,83}(5/3)$.
Here the starting value is 
$3^{4}\cdot 83^{45}$,
a number having $89$ digits. After $2230$ iterations
our algorithm reached the pair $(100,0)$
and hence $$n_{3,83}(5/3) = \lceil3^{100}\cdot3/5\rceil =
309226512439206798621876677859372763621264513201,$$
a number having $48$ digits.
Table \ref{Table-npqalpha} provides some data for further cases and makes manifestly
clear that handling the problem directly - so, 
without working on the exponents only - would be infeasible.

\renewcommand{\arraystretch}{1.3}
\begin{table}[htp] 
\resizebox{\columnwidth}{!}{
\begin{tabular}{ | r | r | c | r | r | r | r | r | r | r | r |}
\hline
$p$ & $q$ & $\alpha$ & $a_\textrm{start}$ & $b_\textrm{start}$  & $a_\textrm{end}$ & $b_\textrm{end}$ & iterations & 
$\ell(n_\textrm{start})$ & $\ell(n_0)$ & $n_{p,q}(\alpha)$\phantom{01234} \\ \hline
$103$ & $223$ & $2$  & $74\,796$ & $6$ & $6$ &  $6\,582$ & $2\,372\,823\,751$ & $150\,567$ & $15\,469$ & $\lceil103^{6}\cdot223^{6582}/2\rceil$ %$11\,905\,643\,330\,881$
\\ \hline
$197$ & $419$ & $2$   & $11\,313$ & $7$ & $7$ & $812$ &  $55\,700\,234$  &  $25\,976$ & $2\,146$ &   $\lceil197^{7}\cdot419^{812}/2\rceil$
\\ \hline 
$13$ & $89$ & $7/4$ &    $44\,875 $ & $4$ & $6$ &  $3\,266$ & $566\,222\,486$ & $49\,997$ & $6\,374$ & $ \lceil 13^{6} \cdot 89^{3266} \cdot4/7\rceil$
\\ \hline 
$137$ & $311$ & $7/4$ & $5\,174$ & $6$ &  $6$ &   $ 1\,410$ & $10\,339\,605$ & $11\,071$ & $3\,528$ & $\lceil 137^{6}\cdot311^{1410} \cdot4/7\rceil$
\\ 
\hline
$89$  & $479$ & $5/3$ & $4\,415$ & $8$ & $10$ & $290$ &  $7\,066\,379$ & $8\,628$ & $797$ &$\lceil 89^{10}\cdot479^{290} \cdot3/5\rceil$
\\ \hline
$293$  & $491$ & $5/3$ & $5\,221$ & $11$ & $11 $ & $55$ &  $12\,553\,989$ & $12\,910$ & $176$ &$\lceil293^{11}\cdot491^{55}\cdot3/5\rceil$
\\ \hline
$79$ & $293$ & $3/2$ & $5\,170$ & $10$ &  $12$ &   $296$ & $10\,276\,001$ & $9\,836$ & $753$ & $\lceil 79^{12}\cdot293^{296} \cdot2/3\rceil$
\\ \hline
$313$ & $487$ & $3/2$ & $14$ & $9\,749$ &  $880$ &   $12$ & $8\,056\,800$ & $26\,236$ & $2\,229$ & $\lceil 313^{880} \cdot487^{12} \cdot2/3\rceil$
\\ \hline
$101$ & $293$ & $4/3$ & $16$ & $14\,616$ &  $3\,412$ &   $12$ & $40\,708\,783$ & $36\,088$ & $6\,869$ & $\lceil 101^{3412} \cdot293^{12} \cdot3/4\rceil$
\\ \hline
$167$ & $367$ & $4/3$ & $29\,063$ & $13$ &  $14$ &   $6\,792$ & $339\,708\,048$ & $64\,633$ & $17\,451$ & $\lceil 167^{14} \cdot367^{6792} \cdot3/4\rceil$
\\ \hline
\end{tabular}} 

\caption{\small{Some data obtained during the computations.
In this table $n_0=p^{a_\textrm{end}}\cdot q^{b_\textrm{end}}$
and  $\ell(n)$ is the number of decimals of $n$.}}
\label{Table-npqalpha}
\end{table} 
We will now give a detailed description of the
determination of $n_{13,89}(7/4)$.
\begin{example}
\label{example-13-89}
To explain why some cases are harder than others, we consider $p=13$ and $q=89$.
The continued fractions convergents of $\theta_{13,89}$ have numerators and denominators respectively
equal to:
\begin{alignat*}{18}
0, & \hskip 0.25cm & 1, & \hskip 0.25cm &  1, & \hskip 0.25cm & 3, & \hskip 0.25cm &   4,  & \hskip 0.25cm & 25643,  & \hskip 0.25cm & 179505, & \hskip 0.25cm & 205148,  
\dotsc
\\
1,& \hskip 0.25cm & 1, & \hskip 0.25cm & 2, & \hskip 0.25cm & 5, & \hskip 0.25cm & 7, & \hskip 0.25cm & 44875, & \hskip 0.25cm & 314132, & \hskip 0.25cm & 359007, 
\dotsc
\end{alignat*}
Hence $r_0/s_0= 0$,  $r_1/s_1= 1$, $r_2/s_2= 1/2$, $r_3/s_3= 3/5$, $r_4/s_4= 4/7$, $r_5/s_5= 25643/44875$ and so on.
If $\alpha=2$, we obtain  $\lceil \frac{\log 89}{\log2} \rceil = 7$ and hence we have  $M=3$,
$s_{M+1}=7$, $r_{M+1}=4$, $s_{M}=5$, $r_{M}=3$ in \eqref{nstart-def}. 
In this case the algorithm starts 
with the information that 
$n_{13,89}(2) \le 
13^5\cdot 89^4 = 23\,295\,754\,887\,613$,
a number having $14$ digits. 
On the other hand  $\lceil \frac{\log 89}{\log(7/4)} \rceil = 9$, $M=4$,
$s_{M+1}=44875$, $r_{M+1}=25643$, $s_{M}=7$, $r_{M}=4$. 
So the algorithm starts just  
with the information that 
$n_{13,89}(7/4) \le 
13^{44875} \cdot 89^4$,
a number having $49\,997$ digits\dots
The huge difference between the height $h(r_4/s_4)$ and 
$h(r_5/s_5)$ is responsible for $n_{13,89}(7/4)$ being much harder to compute than $n_{13,89}(2)$.

The situation for $p=2$ and $q=8191$ is similar:
in this case
$(r_2, s_2)=(1,13)$, $(r_3,s_3)=(5677,73800)$,
$\lceil \frac{\log 8191}{\log2} \rceil = 13$
and $\lceil \frac{\log 8191}{\log (5/3)} \rceil = 23$.
Hence in order to establish that $n_{2,8191}(5/3) = 
\lceil 2^{12} \cdot 8191^{1493}\cdot 3/5 \rceil$,  
the algorithm starts with $n_{\textrm{start}} = 2^{73800}\cdot 8191$, a number having $22\,220$ digits, 
while, to obtain $n_{2,8191}(2) = 1$, it is sufficient
to work with the much smaller starting number $2^{12}\cdot8191 = 33\,550\,336$.
\end{example}

\subsection{Extremal behavior of the gaps}
\label{sect:extremal}
We now show how to analyze the behavior of $D_1$ and $D_2$
implicitly defined in the inequalities
\begin{equation}
\label{D1-D2-def}
\frac{n_i}{(\log n_i)^{D_1}}<n_{i+1}-n_i< \frac{n_i}{(\log n_i)^{D_2}},
\end{equation}
where $n_i \ge 3$.
We further define $\rho_i=\rho_i(p,q)$ implicitly by 
\begin{equation}
\label{rhoi-def}
n_{i+1}-n_i
=
\frac{n_i}{(\log n_i)^{\rho_i}},
\quad i.e.,\quad 
\rho_i = -\frac{\log(\frac{n_{i+1}}{n_i}-1)}{\log \log n_i}
=
\frac{\log(\frac{n_{i}}{n_{i+1}-n_i})}{\log \log n_i},
\end{equation}
so that,  if $n_i \ge 3$, then $D_1 = D_1(p,q) = \max_i \rho_i(p,q)$ and 
$D_2 = D_2(p,q) = \min_i \rho_i(p,q)$.
We note that if 
$n_{j+1}-n_j > n_j$ 
for some $j$ (that is the Bertrand's Postulate property does not hold for $n_j$),  
then  $\rho_j(p,q)$ is negative and hence $D_2(p,q) < 0$.
As this occurs for at most 
finitely many $j$, we can allow ourselves to disregard these outliers.
Hence we will evaluate $\rho_i(p,q)$ only for 
$n_i\ge n_{p,q}(2)$.

We remark that the smallest $D_2$ is zero and this value is reached when $n_{i+1} = 2n_i$;
moreover, the largest  $D_1$ are usually obtained with small powers of primes  and,
in fact, the maximal $D_1$ we got is for $n_i = 3, n_{i+1} = 4$, and is 
$\log3/\log \log 3\approx 11.681421\dotsc$.

We also remark that $D_j = C_j - \log(C_{j+2})/\log \log n_i$, $j=1,2$,
where  $C_1,C_2,C_3, C_4$ are
defined in Theorem \ref{maintijdeman}.
Recalling $C_3 = (\log p)^{C_1}$ and $C_4 = 8q$, this means that
\begin{equation}
\label{C-D-relations}
C_1(p,q) = D_1(p,q) \frac{\log \log n_i}{\log \log n_i - \log\log p}
\quad
\textrm{and}
\quad
C_2(p,q) =  D_2(p,q) + \frac{\log (8q)}{\log \log n_i},
\end{equation}
where, in the first case, we also have to assume $n_i \ne p$.

As we have just explained, in order to have meaningful results we need to work with $\rho_i \ge 0$,
which is equivalent with $n_{i+1} \le 2n_i$. 
We will also require that $n_i := p^a \cdot q^b  \ge 3$ and 
$$
n_i  \ge N:= \exp(\log p\cdot \log q)\,\,(=p^{\log q}\,=q^{\log p}).
$$
Since sometimes we do not know the value of $n_{p,q}(2)$,
the best strategy is then to start working with $n_\textrm{start}$ as identified 
in Section \ref{nstart-gen} and use a modified form of
the left neighbor search (explained in Section \ref{left-neighbors-gen}) to generate
$n_{i-1}$, until  we have produced  $L$ left neighbors, or 
we have reached  $n_{p,q}(2)$.
To do this, we need to slightly modify the search procedure in Section \ref{left-neighbors-gen},
since in this case the issue is 
to determine $n_{i-1}$, rather than to get as close as 
possible to $n_{p,q}(2)$.
Hence in step \ref{cfrac-search}  of Section \ref{left-neighbors-gen} 
we will always choose the value of $x_{\textrm{max}}$ 
and in steps \ref{trivial-first-search}-\ref{trivial-second-search} of the same section
we will search for the maximal value of the form $x := a\log p + b \log q$
in the ranges there defined,
where $a,b\in {\mathbb Z}$ run in the intervals
$0\le b \le \lfloor\log n_{i}/\log q\rfloor$ and
$0\le a \le \lfloor(\log n_{i} - b \log q)/\log p\rfloor$. These
changes to the search procedures increase their computational cost;
unfortunately, there is no way to work differently since in this case
we cannot halt the procedure as soon as we have found an $S$-unit less than $n_i$
in the prescribed interval, but we need to be sure that such a point is the maximal 
$S$-unit less than $n_i$, or, in other words, that such a point is the left neighbor $n_{i-1}$
of $n_i$.

If we have reached $n_{p,q}(2)$ without having generated $L$ left neighbors,
we start to generate right neighbors from $n_\textrm{start}$ until 
we have obtained a total number of $L$ neighbors.
For the generation of the right neighbor $n_{i+1}$ of $n_i$, we used
the continued fraction approach described in \cite[Theorem 2.2]{BDH-2014}.
The theoretical justification for the neighbor search procedure is provided
in the Appendix.

In each step of the previously described algorithm, 
as soon as we have determined one of the neighbors of $n_i$,
we can evaluate $\rho_i$, defined in \eqref{rhoi-def}, the constants defined in \eqref{D1-D2-def}
and   \eqref{C-D-relations}. 

The same remarks we made in Section \ref{sect:comp_npq_alpha}
about using the exponents only in the computations apply here 
as well with a single exception: in evaluating $\log(n_{i+1}/ n_i -1)$ 
in \eqref{rhoi-def}, we are in fact forced  either to generate both 
$n_i$ and $n_{i+1}$, or to use 
\[
\log(n_{i+1}/{n_i} -1) 
= 
\log \bigl( \exp(u \log p + v \log q ) - 1 \bigr),
\]
in which we assume $u \log p + v \log q$, 
where $u,v\in \Z$ are such that $n_{i+1}/ n_i = p^{u}\cdot q^{v}$, to be known.
This is one of the most computationally costly steps of the whole procedure,
but unfortunately there is no other way to compute $\rho_i$.

In this way we were able to collect, for every $2\le p<q$,
with $p$ and $q$ both primes, all the
values of $C_1,C_2,C_3,C_4,D_1,D_2$ 
and of their averaged values (defined in the next section).
This is a heavy computation
and the largest $P$ we were able to work with was $P=10^4$.
The data in Table \ref{Tablec1c2mu} were obtained 
as a part of computation having an accuracy of
$19$ decimal digits (but the results are here truncated at $10$ digits) 
performed for every $2\le p<q< 10^4$,
and using \genneighbor\ neighbors, 
a computation which required about 
$2$ days and $18$ hours on the Dell Optiplex machine mentioned before.

\renewcommand{\arraystretch}{1.3}
\begin{table}[htp] 
\resizebox{\columnwidth}{!}{
\begin{tabular}{ | r | r | l | l | l | }
\hline
$p$ & $q$ & \phantom{01}$D_{1}(p,q)$ & \phantom{01}$D_{2}(p,q)$ & \phantom{01}$ \mu(p,q;k)$ \\ \hline
$2$ & $3$ & $11.6814212141$ & $0.6261567724$ & $1.0941489589$ \\ \hline
$2$ & $5$ &  $4.2441792885$ & $0.3589187149$ & $0.9474929648$ \\ \hline
$2$ & $7$ &  $2.9229728309$ & $0.1201698368$ & $0.8284546634$ \\ \hline
$2$ & $11$ & $1.8178288541$ & $0.2942087461$ & $0.9858711527$ \\ \hline
$2$ & $13$ & $1.7851594152$ & $0.2088055868$ & $0.7625506652$ \\ \hline
$2$ & $17$ & $2.7188068070$ & $0$            & $0.8044481097$ \\ \hline
$2$ & $19$ & $2.0907970943$ & $0$            & $0.7925235409$ \\ \hline
$2$ & $23$ & $1.5215514442$ & $0.2221675922$ & $0.6835338899$ \\ \hline
$2$ & $29$ & $1.8685980225$ & $0.0635501877$ & $0.6305104841$ \\ \hline
$2$ & $31$ & $2.7834367088$ & $0.0149561107$ & $0.8398947566$ \\ \hline
$2$ & $37$ & $1.4934915642$ & $0$            & $0.8004576906$ \\ \hline
$2$ & $41$ & $1.3988710027$ & $0$            & $0.7933541465$ \\ \hline
$2$ & $43$ & $1.3718840668$ & $0$            & $0.7791951977$ \\ \hline
$2$ & $47$ & $1.4407350298$ & $0$            & $0.6686016472$ \\ \hline
$2$ & $53$ & $1.4055386972$ & $0$            & $0.6391686515$ \\ \hline
$2$ & $59$ & $1.7560379970$ & $0.0474250925$ & $0.7734676479$ \\ \hline
$2$ & $61$ & $2.1308613919$ & $0.0239511431$ & $0.7720445935$ \\ \hline
$2$ & $67$ & $2.1471869534$ & $0.0460488967$ & $0.6122509704$ \\ \hline
$2$ & $71$ & $1.5526948207$ & $0.1074029002$ & $0.7469833108$ \\ \hline
$2$ & $73$ & $1.3763643264$ & $0.1382120440$ & $0.7668223147$ \\ \hline
$2$ & $79$ & $1.0720303269$ & $0.2320129747$ & $0.7474188012$ \\ \hline
$2$ & $83$ & $2.3638631955$ & $0.2964674293$ & $0.8588722462$ \\ \hline
$2$ & $89$ & $1.5378684405$ & $0.2093917801$ & $0.7841929598$ \\ \hline
$2$ & $97$ & $2.5893808896$ & $0.2590800786$ & $0.6246448440$ \\ \hline
\end{tabular}
\begin{tabular}{ | r | r | l | l | l | }
\hline
$p$ & $q$ & \phantom{01}$D_{1}(p,q)$ & \phantom{01}$D_{2}(p,q)$ & \phantom{01}$ \mu(p,q;k)$ \\ \hline
$3$ & $5$  & $2.1605296213$ & $0.1323896610$ & $0.8603600118$ \\ \hline
$3$ & $7$  & $1.1710669001$ & $0.1383426171$ & $0.7649429172$ \\ \hline
$3$ & $11$ & $1.4012089272$ & $0.2086747873$ & $0.7566780823$ \\ \hline
$3$ & $13$ & $2.6410842719$ & $0.1285354272$ & $0.7479745611$ \\ \hline
$3$ & $17$ & $1.7461871466$ & $0.0795647361$ & $0.7255147576$ \\ \hline
$3$ & $19$ & $1.4225074511$ & $0.2039287543$ & $0.7949565378$ \\ \hline
$3$ & $23$ & $1.2065704232$ & $0.0656360747$ & $0.7683872665$ \\ \hline
$3$ & $29$ & $0.9394181270$ & $0.0151595599$ & $0.7455240234$ \\ \hline
$3$ & $31$ & $1.5149281514$ & $0.0071479409$ & $0.6700758602$ \\ \hline
$3$ & $37$ & $1.4807635709$ & $0.2296184840$ & $0.6676913346$ \\ \hline
$3$ & $41$ & $1.0229032857$ & $0.0138973896$ & $0.7159407244$ \\ \hline
$3$ & $43$ & $0.9916232421$ & $0.0692590381$ & $0.7018910831$ \\ \hline
$3$ & $47$ & $2.2545776040$ & $0.0535541917$ & $0.7078101517$ \\ \hline
$3$ & $53$ & $0.9771181314$ & $0.0200128231$ & $0.7076989903$ \\ \hline
$3$ & $59$ & $1.5682950986$ & $0.2189707656$ & $0.8062334056$ \\ \hline
$3$ & $61$ & $1.2885326430$ & $0.1479792903$ & $0.8006209830$ \\ \hline
$3$ & $67$ & $0.9915469462$ & $0.1316451029$ & $0.7048711907$ \\ \hline
$3$ & $71$ & $1.3780698058$ & $0.0871486016$ & $0.7826861503$ \\ \hline
$3$ & $73$ & $0.9560032551$ & $0.0070696424$ & $0.6675886273$ \\ \hline
$3$ & $79$ & $1.2434176789$ & $0.0090970147$ & $0.7737402410$ \\ \hline
$3$ & $83$ & $1.3056800194$ & $0.0042329888$ & $0.6944828831$ \\ \hline
$3$ & $89$ & $1.4320941282$ & $0.0417272844$ & $0.7999629578$ \\ \hline
$3$ & $97$ & $1.2268141005$ & $0.1029715740$ & $0.6766541024$ \\ \hline
\phantom{} & \phantom{} & \phantom{}  & \phantom{}  & \phantom{}\\ \hline
\end{tabular}
\begin{tabular}{ | r | r | l | l | l | }
\hline
$p$ & $q$ & \phantom{01}$D_{1}(p,q)$ & \phantom{01}$D_{2}(p,q)$ & \phantom{01}$ \mu(p,q;k)$ \\ \hline
$5$ & $7$ &  $1.1384784153$ & $0.1109894596$ & $0.7781281263$ \\ \hline
$5$ & $11$ & $1.1887010049$ & $0.0014698384$ & $0.6590355370$ \\ \hline
$5$ & $13$ & $1.3211241477$ & $0.0310425989$ & $0.7757032047$ \\ \hline
$5$ & $17$ & $1.1035830903$ & $0.1797045904$ & $0.6901667574$ \\ \hline
$5$ & $19$ & $1.1398389151$ & $0.1578288483$ & $0.6723537325$ \\ \hline
$5$ & $23$ & $0.9315775235$ & $0.0005087729$ & $0.6645526054$ \\ \hline
$5$ & $29$ & $1.0391485562$ & $0.0813597681$ & $0.6529045038$ \\ \hline
$5$ & $31$ & $1.2305570961$ & $0.1189302077$ & $0.5509941851$ \\ \hline
$5$ & $37$ & $1.8067753639$ & $0.1479331658$ & $0.7585371271$ \\ \hline
$5$ & $41$ & $1.2902483645$ & $0.0690378980$ & $0.5759361463$ \\ \hline
$5$ & $43$ & $1.6660279429$ & $0.0636126574$ & $0.7032472651$ \\ \hline
$5$ & $47$ & $0.9651734788$ & $0.0491064048$ & $0.6754674325$ \\ \hline
$5$ & $53$ & $1.2873497260$ & $0.0322254510$ & $0.4291499629$ \\ \hline
$5$ & $59$ & $1.3174229255$ & $0.0330207317$ & $0.6428076307$ \\ \hline
$5$ & $61$ & $1.0939162195$ & $0.1045174998$ & $0.6276138163$ \\ \hline
$5$ & $67$ & $1.0048592858$ & $0.0534692090$ & $0.6291197333$ \\ \hline
$5$ & $71$ & $1.0184013877$ & $0.1011594660$ & $0.6371897583$ \\ \hline
$5$ & $73$ & $2.1481873755$ & $0.0470304122$ & $0.5549093481$ \\ \hline
$5$ & $79$ & $1.4575444151$ & $0.0011909162$ & $0.5340714572$ \\ \hline
$5$ & $83$ & $1.2343218159$ & $0.1386930533$ & $0.5927933668$ \\ \hline
$5$ & $89$ & $0.9266169776$ & $0.0781941449$ & $0.6043041547$ \\ \hline
$5$ & $97$ & $1.0334458810$ & $0.0680491434$ & $0.5810298221$ \\ \hline
\phantom{} & \phantom{} & \phantom{}  & \phantom{}  & \phantom{}\\ \hline
\phantom{} & \phantom{} & \phantom{}  & \phantom{}  & \phantom{}\\ \hline
\end{tabular}
}
\caption{\small{Computed values of $D_1(p,q),D_2(p,q), \mu(p,q;k)$ with $2\le p \le 5$, 
$p < q \le 97$ and having generated \genneighbor\ neighbors
for each case.}}
\label{Tablec1c2mu} 
\end{table}

\subsection{Average behavior of the gaps} 
Define 
\[
\mu(p,q;k):= \frac{1}{k^+}
\sum_{\substack{i=1\\\rho_i(p,q)\ge 0}}^k \rho_i(p,q),\quad{with} 
\ k^+:=\#\{1\le i\le k:\rho_i(p,q)\ge 0\}.
\]
\begin{lemma}
Let $p,q$ be fixed and $\epsilon>0$. 
Recall that $C_1=2\cdot 10^9\log p\cdot \log q$.
There exists an integer $k_{p,q}(\epsilon)$ such
that $$C_1^{-1}-\epsilon <\mu(p,q;k)<C_1\quad\textrm{for~every}\ k\ge k_{p,q}(\epsilon).$$
\end{lemma}
\begin{Proof}
By Theorem \ref{maintijdeman} we have
$$C_1^{-1}-\frac{\log C_4}{\log \log n_i}< \rho_i< C_1-\frac{\log C_3}{\log \log n_i}.$$
The proof follows from these two inequalities, $C_3>0$ and the observation that $(\log \log n_i)^{-1}$ tends to zero.
\end{Proof}
We cannot answer the following natural question.
\begin{question}
Does $\lim_{k\rightarrow \infty}\mu(p,q;k)$ exist?
\end{question}
The numerical work presented here suggests that the answer is \textbf{yes}. In this case we write $\mu(p,q)$ for the limit.
\begin{remark}
Instead of
$\mu(p,q;k)$ one can consider $\frac{1}{k}\sum_{i=1}^k\rho_i(p,q)$, which
has also limit $\mu(p,q)$, if the limit exists. However, our preference is 
to work with $\mu(p,q;k)$, as numerically it seems to behave more regularly.
\end{remark}

\smallskip
The programs and the results here described are available at the address:\newline
\url{www.math.unipd.it/~languasc/Sunits.html}.

\medskip
\noindent {\tt Acknowledgement.} 
We thank the authors of \cite{BDH-2014} for helpful correspondence.
\par A substantial part of the work described
here, was carried out during the research stay of Togb\'e at the
Max-Planck Institute for Mathematics in June-July 2023. A short
visit of Luca also proved to be very helpful.
Luca, Moree and Togb\'e thank the MPIM and its
staff for the great working environment, the hospitality and the support.
\par The discriminator values reported on were kindly computed by Matteo Ferrari. The integers $n_{p,q}(\alpha)$ were computed on a machine
of the cluster located at the Dipartimento di Matematica ``Tullio Levi-Civita'' of 
the University of Padova, see \url{https://hpc.math.unipd.it}. 
Languasco is grateful for having had such computing facilities 
at his disposal.

\appendix 
\section{On efficient left and right neighbors searches}
\label{App-searches-thms}
Here we establish two theorems on
determining the neighbors of $n_i$.
We will just make use of \cite[Lemma 3.1]{BDH-2014}, which characterizes
the best approximants of a real number $\theta$ in terms of the principal and intermediate
convergents generated by its continued fraction. 
In line with \cite{BDH-2014} we will use the notation $s_i$ instead of $n_i$. 
The sought for neighbors are then $s_{i-1}$ and $s_{i+1}$.
\begin{theorem}
[Left neighbor search]
\label{predecessor-thm}
Suppose $k\ge 1$ and that we are given $s_{k} = p^{c_k} \cdot q^{d_k}$,
$(c_k,d_k) \in{\mathbb Z}_{\ge 0}^2$, $(c_k,d_k) \ne (0,0)$.
Then we can compute its left neighbor $s_{k-1}$ in the following way:
\begin{enumerate}[leftmargin=*, label={\roman*)}]
\item Let $u_1/v_1>0$, $(u_1,v_1)\in {\mathbb Z}_{\ge 0}^2$, be the upper convergent of $\theta_{p,q}$ 
with maximal numerator for which $u_1 \le d_{k}$ holds.
\item Let $u_2/v_2>0$, $(u_2,v_2)\in{\mathbb Z}_{\ge 0}^2$, be the lower convergent of $\theta_{p,q}$ 
with maximal denominator for which $v_2 \le c_{k}$ holds.
\item Put $x: = \vert v_1 \log p - u_1 \log q \vert - \vert v_2 \log p - u_2 \log q \vert$
and
\[
c_{k-1} = 
\begin{cases}[1]
c_k + v_1 & \textrm{if}\ \ x<0,\\
c_k - v_2 & \textrm{if}\ \ x>0, 
\end{cases}
\qquad
d_{k-1} = 
\begin{cases}[1]
d_k - u_1 & \textrm{if}\ \ x<0,\\
d_k + u_2 & \textrm{if}\ \ x>0.
\end{cases}
\]
\end{enumerate}
Then we have $s_{k-1} = p^{c_{k-1}}\cdot q^{d_{k-1}}$.
\end{theorem}  

\begin{Proof}
Let  $x<0$. Then
we have $c_{k-1} = c_k + v_1 $, $d_{k-1} = d_k - u_1$.
Since $u_1/v_1$ is an upper convergent of $\theta_{p,q}= \log p/\log q$,
we have $v_1 \log p  -u_1 \log q < 0$, and we define 
$S :=  s_k p^{v_1} \cdot q^{-u_1} < s_k$.
Now we prove that $S= s_{k-1}$.

By contradiction, if $S \ne s_{k-1}$ 
there exists $(a,b)\in{\mathbb Z}_{\ge 0}^2$,  $(a,b) \ne (0,0)$, such that  
$S = p^{c_k+v_1} \cdot q^{d_k-u_1} < p^a \cdot q^b < p^{c_k} \cdot q^{d_k} = s_k$.
This is equivalent to 
\begin{equation}
\label{fund1}
0
< 
\frac{d_k-b}{\theta_{p,q}}  - (a - c_k)  
< 
\frac{ u_1}{\theta_{p,q}}   - v_1. 
\end{equation}

If $b=d_k$,  we must have $a<c_k$. Moreover,  we will obtain that $u_1/(v_1+c_k-a)$ is an upper
convergent of $\theta_{p,q}$ having the same numerator of $u_1/v_1$.
This implies $a=c_k$, which is impossible because $p^a \cdot q^b < s_k$ by definition.

So we have that $b\ne d_k$. Now, thanks to  \cite[Lemma 3.1]{BDH-2014}, 
from \eqref{fund1} we obtain $d_k-b > u_1$.
If $b > d_k$, it follows that $d_k  > b+u_1> d_k +u_1$ and hence $u_1<0$,
which is a contradiction.
Assume now $0 \le b < d_k$. Then $\overline{u} := d_k-b > 0$
and $\overline{u} > u_1$. Again using \eqref{fund1}, we have  
$\overline{v} := a-c_k > (\overline{u} - u_1)/\theta_{p,q} + v_1 > 0$
and $\overline{u} /\overline{v} $ is an upper convergent for $\theta_{p,q}$.
But this is a contradiction, since $u_1$ is maximal between the numerators $\le d_k$.
This proves that $s_{k-1} = p^{c_k+v_1} \cdot q^{d_k-u_1}$.
The case  $x>0$ can be proved analogously.
\end{Proof}
 
In Theorem \ref{predecessor-thm}  we used the same notation used in Theorem 2.2 in \cite{BDH-2014};
we think that this would help the reader to spot their differences more easily.
Remark that the choices of the convergents of $\theta_{p,q}$  in this theorem
are precisely the ones in Section \ref{left-neighbors-gen}, point \ref{cfrac-search},
and hence this $s_{k-1}$ corresponds to the choice of $x_{\textrm{max}}$ there.

We also include an alternative, and shorter, proof of Theorem 2.2 in \cite{BDH-2014}.
\begin{theorem}[Right neighbor search]
\label{successor-thm-complete}
Suppose $k\ge 0$ and that we are given $s_{k} = p^{c_k} \cdot q^{d_k}$,
$(c_k,d_k) \in{\mathbb Z}_{\ge 0}^2$.
Then we can compute its right neighbor $s_{k+1}$ in the following way:
\begin{enumerate}[leftmargin=*, label={\roman*)}]
\item Let $u_1/v_1>0$, $(u_1,v_1)\in{\mathbb Z}_{\ge 0}^2$, be the upper convergent of $\theta_{p,q}$ 
with maximal denominator for which $v_1 \le c_{k}$ holds.
\item Let $u_2/v_2>0$, $(u_2,v_2) \in{\mathbb Z}_{\ge 0}^2$, be the lower convergent of $\theta_{p,q}$ 
with maximal numerator for which $u_2 \le d_{k}$ holds.
\item Put $x: = \vert v_1 \log p - u_1 \log q \vert - \vert v_2 \log p - u_2 \log q \vert$
and 
\begin{equation*}
c_{k+1} = 
\begin{cases}[1]
c_k - v_1 & \textrm{if}\ x<0,\\
c_k + v_2 & \textrm{if}\ x>0, 
\end{cases}
\qquad
d_{k+1} = 
\begin{cases}[1]
d_k + u_1 & \textrm{if}\ x<0,\\
d_k - u_2 & \textrm{if}\ x>0.
\end{cases}
\end{equation*} 
\end{enumerate}
Then we have $s_{k+1} = p^{c_{k+1}}\cdot q^{d_{k+1}}$.
\end{theorem}  

\begin{Proof}
Assume  $x<0$. Then
$ \vert v_1 \log p - u_1 \log q \vert < \vert v_2 \log p - u_2 \log q \vert$
and we choose $c_{k+1} = c_k - v_1 $, $d_{k+1} = d_k + u_1$.
Since $u_1/v_1$ is an upper convergent of $\theta_{p,q}= \log p/ \log q$,
we have $u_1 \log q  - v_1 \log p >  0$, and we define 
$S :=  s_k p^{-v_1} \cdot q^{u_1} > s_k$.
Now we prove that $S= s_{k+1}$. 

By contradiction, if $S \ne s_{k+1}$ 
there exists $a,b\ge 0$,  $(a,b) \ne (0,0)$, such that  
$ p^{c_k-v_1} \cdot q^{d_k+ u_1}=S  >  p^a \cdot q^b > s_k = p^{c_k} \cdot q^{d_k}$.
This is equivalent to 
\begin{equation}
\label{fund3-appendix}
0<
 (b-d_k) - (c_k -a) \theta_{p,q}
<
 u_1   - v_1 \theta_{p,q} .
\end{equation}

If $a=c_k$,  we must have $b > d_k$. Moreover,  we will obtain that $(d_k-b + u_1)/v_1$ is an upper
convergent of $\theta_{p,q}$ having the same denominator of $u_1/v_1$.
This implies $b=d_k$, which is impossible because $p^a \cdot q^b > s_k$ by definition. 

So we have that $a\ne c_k$. 
Now, thanks to  \cite[Lemma 3.1]{BDH-2014}, from \eqref{fund3-appendix} 
we obtain $c_k  > a+v_1$.
If $a > c_k$, we obtain $c_k  > a+v_1> c_k +v_1$ and hence $v_1<0$,
which is a contradiction.
Assume now $0 \le a < c_k$. Then $\overline{v} := c_k-a > 0$
and $\overline{v} > v_1$. 
Again using \eqref{fund3-appendix}, we have   $\overline{u}:= b-d_k > \overline{v} \theta_{p,q}>0$ 
and $\overline{u} /\overline{v} $
is an upper convergent for $\theta_{p,q}$.
But this is a contradiction, since $v_1$ is maximal between the denominators $\le c_k$.
This proves that $s_{k+1} = p^{c_k-v_1}  \cdot q^{d_k+u_1}$.
The case $x>0$ can be proved analogously.
 \end{Proof}

\end{document}